\documentstyle[12pt]{article}

\def\qed{$\Box$}
\def\colon{\;{:}\;}

\def\int{\intop\limits}
\def\eqref#1{\rm (\ignorespaces\ref{#1}\unskip)}
\let\cases\relax
\def\cases{\left\{\begin{array}{cl}}
\def\endcases{\end{array}\right.}
\def\text#1{\quad\mbox{#1}\quad}
\def\tfrac#1#2{\textstyle\frac{#1}{#2}}
\def\subjclass#1{{\def\thefootnote{}\footnotetext{1991 {\it Mathematics
      Subject Classification.}\ #1}}}

\def\t{\mbox{\boldmath$\tau$}}
\let\rho\varrho
\let\phi\varphi
\def\height{\mathop{\rm lh}\nolimits}

\font\twelvemsb=msbm10 scaled 1200
\font\tenmsb=msbm10
\font\eightmsb=msbm8
\font\sevenmsb=msbm7
\font\sixmsb=msbm6
\font\fivemsb=msbm5
\font\twelveeusm=eusm10 scaled 1200

\font\twelveeufb=eufb10 scaled 1200

\makeatletter
  \ifcase\@ptsize                      %
    \newfam\msbfam                     %
    \textfont\msbfam=\tenmsb           
    \scriptfont\msbfam=\sevenmsb       
    \scriptscriptfont\msbfam=\fivemsb  %
    \def\eusm#1{\mbox{\teneusm #1}}    %
    \def\frak#1{\mbox{\teneufb #1}}    %
  \or
  \or
    \newfam\msbfam                     %
    \textfont\msbfam=\twelvemsb        %
    \scriptfont\msbfam=\eightmsb       
    \scriptscriptfont\msbfam=\sixmsb   
    \def\eusm#1{\mbox{\twelveeusm #1}} %
    \def\frak#1{\mbox{\twelveeufb #1}} %
  \fi
\makeatother

\def\Bbb#1{{\fam\msbfam\relax#1}}
\def\B{\Bbb B\,}
\def\D{\Bbb D\,}
\def\F{\Bbb F\,}

\def\N{\Bbb N}

\def\S{\Bbb S}

\makeatletter
  \def\@begintheorem#1#2{\it 
    \trivlist \item[\hskip \labelsep{\bf #1\ #2.}]}
  \def\@opargbegintheorem#1#2#3{\it 
    \trivlist \item[\hskip \labelsep{\bf #1\ #2\ (#3).}]}
\makeatother

\newtheorem{corollary}{Corollary}       
\newtheorem{proposition}[corollary]{Proposition}
\newtheorem{theorem}[corollary]{Theorem}
\newtheorem{lemma}[corollary]{Lemma}

\newenvironment{definition}{\ifvmode\else\newline\fi{\bf Definition.}}{}

\def\remarkname{Remark}
\newenvironment{remark}{\ifvmode\else\newline\fi{\it \remarkname.}}{}

\newcounter{problem}
\def\theproblem{{\rm \arabic{problem}}}

\newenvironment{problem}{\refstepcounter{problem}\ifvmode\else\newline\fi{\it
    Problem \theproblem.}}{}

\def\proofname{Proof}
\def\pf{{\sc \proofname:\ }}

\begin{document}

\title{Vector--valued Walsh--Paley martingales and geometry of Banach
  spaces}
\author{J\"org Wenzel}

\maketitle

\begin{abstract}
The concept of Rademacher type $p$ ($1\leq p\leq2$) plays an important
role in the local theory of Banach spaces. 
In \cite{mas88} Mascioni
considers a weakening of this concept and shows that for a Banach
space $X$ weak Rademacher type $p$ implies Rademacher type $r$ for all
$r<p$.

As with Rademacher type $p$ and weak Rademacher type $p$, we introduce
the concept of Haar type $p$ and weak Haar type $p$ by replacing the
Rademacher functions by the Haar functions in the respective
definitions. 
We show that weak Haar type $p$ implies Haar type $r$ for all $r<p$. 
This solves a problem left open by Pisier \cite{pis75}.

The method is to compare Haar type ideal norms related to different
index sets.
\end{abstract}


\section{Introduction}
\label{sec:1}

Let $(r_n)$ denote the sequence of Rademacher functions. 
For $1\leq p\leq2$, we say that a Banach space $X$ is of {\em
  Rademacher type\/} $p$ if there exists a constant $c\geq0$ such that
$$ 
  \bigg\|\sum_{k=1}^nx_kr_k\bigg| L_p\bigg\|\leq c\bigg( \sum_{k=1}^n
  \|x_k\|^p\bigg)^{1/p} 
$$ 
for all $n\in\N$ and $x_1,\dots,x_n\in X$.

Furthermore, for $1\leq p<2$, we say that $X$ is of {\em weak Rademacher
  type\/} $p$ if there exists a constant $c\geq0$ such that
$$ 
  \bigg\|\sum_{k=1}^nx_kr_k\bigg| L_2\bigg\|\leq c\,n^{1/p-1/2}
  \bigg( \sum_{k=1}^n \|x_k\|^2\bigg)^{1/2} 
$$ 
for all $n\in\N$ and $x_1,\dots,x_n\in X$.

Kahane's Inequality and the Cauchy--Schwarz Inequality tell us that
Rademacher type $p$ implies  weak Rademacher type $p$. 
Mascioni \cite{mas88} has shown a partial converse.
\begin{theorem}[Mascioni \cite{mas88}]\label{th:1}
  If a Banach space $X$ is of weak Rademacher type $p$ for some $1\leq
  p<2$ then it is of Rademacher type $r$ for all $r<p$.
\end{theorem}

A different concept can be introduced by using sequences of
$X$--valued martingale differences instead of Rademacher sequences.
Following Pisier \cite{pis86}, we say that a Banach space $X$ is of {\em
martingale type\/} $p$ ($1\leq p\leq2$) if there exists a constant
$c\geq0$ such that
$$ 
  \bigg\|\sum_{k=1}^nd_k\bigg| L_p\bigg\|\leq c\bigg( \sum_{k=1}^n
  \|d_k|L_p\|^p\bigg)^{1/p} 
$$ 
for all $n\in\N$ and all $X$--valued martingale difference sequences
$d_1,\dots,d_n$.

Again, for $1\leq p<2$, we say that $X$ is of {\em weak martingale
  type\/} $p$ if there exists a constant $c\geq0$ such that
$$ 
  \bigg\|\sum_{k=1}^nd_k\bigg| L_2\bigg\|\leq c\,n^{1/p-1/2} \bigg(
  \sum_{k=1}^n \|d_k|L_2\|^p\bigg)^{1/p} 
$$ 
for all $n\in\N$ and all $X$--valued martingale difference sequences
$d_1,\dots,d_n$.

In \cite{pis75} Pisier showed a theorem similar to Theorem \ref{th:1}.
\begin{theorem}[Pisier \cite{pis75}]\label{th:2}
  If a Banach space $X$ is of weak martingale type $p$ for some $1\leq
  p<2$ then it is of martingale type $r$ for all $r<p$.
\end{theorem}

Denoting by $\chi_k^{(j)}$ the Haar functions, we see that the
sequence
$$ 
  \bigg( \sum_{j=1}^{2^{k-1}}x_k^{(j)} \chi_k^{(j)} \bigg) \quad
  k=1,\dots,n\quad (x_k^{(j)}) \subseteq X 
$$ 
forms a sequence of $X$--valued martingale differences. 
Restricting to those martingales (usually referred to as dyadic or
Walsh--Paley martingales) in the definition of martingale type, we can
define the apparently weaker concept of Haar type. 
For $1\leq p\leq2$, we say that a Banach space $X$ is of {\em Haar
  type\/} $p$ if there exists a constant $c\geq0$ such that
$$ 
  \bigg\|\sum_{k=1}^n \sum_{j=1}^{2^{k-1}} x_k^{(j)}
  \chi_k^{(j)}\bigg| L_p\bigg\|\leq c\bigg( \sum_{k=1}^n
  \bigg\|\sum_{j=1}^{2^{k-1}}x_k^{(j)} \chi_k^{(j)}\bigg| L_p
  \bigg\|^p\bigg)^{1/p} 
$$ 
for all $n\in\N$ and $(x_k^{(j)})\subseteq X$.

\pagebreak[3]
Once more it was Pisier \cite{pis75} who showed that Haar type $p$ and
martingale type $p$ coincide.
\begin{theorem}[Pisier \cite{pis75}]
  A Banach space $X$ is of martingale type $p$ if and only if it is of
  Haar type~$p$.
\end{theorem}

Of course, the next step is to define weak Haar type. 
For $1\leq p <2$, we say that a Banach space $X$ is of {\em weak Haar
  type\/} $p$ if there exists a constant $c\geq0$ such that
$$ 
  \bigg\|\sum_{k=1}^n \sum_{j=1}^{2^{k-1}} x_k^{(j)} \chi_k^{(j)}
  \bigg| L_2\bigg\| \leq c\,n^{1/p-1/2}\bigg( \sum_{k=1}^n
  \bigg\| \sum_{j=1}^{2^{k-1}} x_k^{(j)} \chi_k^{(j)} \bigg| L_2
  \bigg\|^2 \bigg)^{1/2} 
$$ 
for all $n\in\N$ and $(x_k^{(j)})\subseteq X$.

The main result of this paper is the following companion to
Theorems~\ref{th:1} and~\ref{th:2}.
\begin{theorem}\label{th1}
  If a Banach space $X$ is of weak Haar type $p$ for some $1\leq p<2$
  then it is of Haar type $r$ for all $r<p$.
\end{theorem}

This solves a problem in \cite{pis75}, where Pisier could show the
same conclusion under the stronger assumption that $X$ is of weak
martingale type $p$.

Let us now quickly review the contents of this article section by
section.

In Section \ref{sec:2}, we define the necessary concepts. 
In Section \ref{sec:3}, we formulate and prove the main new ingredient
of the proof of the main theorem. 
In Section \ref{sec:3a} we provide a lemma of combinatorial character
needed to prove a local variant of the main theorem.  
In Section \ref{sec:4}, we show how to derive the main theorem from
the results in Section \ref{sec:3}. 
Finally, in Section \ref{sec:final}, we consider some examples and
formulate some problems. 

I am grateful to Albrecht Pietsch for his numerable hints and useful
remarks, which helped to smooth out the content of this paper.


\section{Definitions}
\label{sec:2}

For $k=1,2,\dots$ and $j=0,\pm1,\pm2,\dots$, we define the {\em Haar
functions} by
$$
  \chi_k^{(j)}(t):=
  \begin{cases}
    +2^{(k-1)/2} & \text{ if $t\in\Delta_k^{(2j-1)}$,}\\
    -2^{(k-1)/2} & \text{ if $t\in\Delta_k^{(2j)}$,}\\
    0 & \text{ otherwise.}
  \end{cases}
$$
Here 
$$
  \Delta_k^{(j)}:=\left[\tfrac{j-1}{2^k},\tfrac j{2^k}\right)
$$
are the {\em dyadic intervals}.

The following facts are obvious consequences of the definition of the
Haar functions
\begin{eqnarray}
  \label{eq2a}
  \chi_k^{(j)}(t-\tfrac1{2^{k-1}}) & = & \chi_k^{(j+1)}(t) \quad \text{ and } 
  \\ 
  \chi_{k+1}^{(j)}(t) & = & \sqrt2\chi_k^{(j)}(2t).
\end{eqnarray}

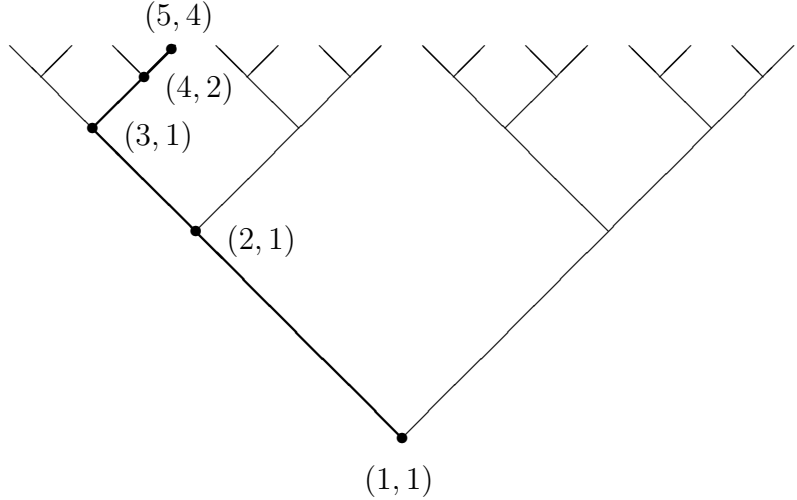
\begin{figure}[ht]
  \begin{center}
    \leavevmode
    \setlength{\unitlength}{0.018in}%
    %
    \begin{picture}(228,130)(96,560)
      \thinlines 
      \put(210,570){\circle*{3}} \put(150,630){\circle*{3}}
      \put(120,660){\circle*{3}} \put(135,675){\circle*{3}}
      \put(143,683){\circle*{3}} 
      \put(210,570){\line( 1, 1){114}} \put(300,660){\line(-1, 1){ 24}} 
      \put(270,630){\line(-1, 1){ 54}} \put(240,660){\line( 1, 1){ 24}}
      \put(210,570){\line(-1, 1){114}} \put(315,675){\line(-1, 1){  9}}
      \put(285,675){\line( 1, 1){  9}} \put(255,675){\line(-1, 1){  9}}
      \put(225,675){\line( 1, 1){  9}} \put(195,675){\line(-1, 1){  9}}
      \put(165,675){\line( 1, 1){  9}} \put(135,675){\line(-1, 1){  9}}
      \put(105,675){\line( 1, 1){  9}} \put(120,660){\line( 1, 1){ 24}}
      \put(150,630){\line( 1, 1){ 54}} \put(180,660){\line(-1, 1){ 24}}
      \thicklines
      \put(210,570){\line(-1, 1){ 60}} \put(150,630){\line(-1, 1){ 30}}
      \put(120,660){\line( 1, 1){ 15}} \put(135,675){\line( 1, 1){  9}}
      \thinlines
      \put(199,555){%
        \makebox(0,0)[lb]{$(1,1)$}}
      \put(159,625){%
        \makebox(0,0)[lb]{$(2,1)$}}
      \put(129,655){%
        \makebox(0,0)[lb]{$(3,1)$}}
      \put(141,669){%
        \makebox(0,0)[lb]{$(4,2)$}}
      \put(135,690){%
        \makebox(0,0)[lb]{$(5,4)$}}
    \end{picture}
  \end{center}
  \caption{The dyadic tree $\D$}
  \label{fig:dyadic tree}
\end{figure}

We denote by 
$$
  \D:=\{(k,j)\colon k=1,2,\dots;j=1,\dots,2^{k-1}\}
$$
the {\em dyadic tree}. 
For a finite subset $\F\subseteq\D$, we consider the orthonormal
system
$$
  \eusm H(\F):=\{\chi_k^{(j)}\colon (k,j)\in\F\}\subseteq L_2[0,1).
$$ 
In particular, we will consider finite dyadic trees
$$
  \D_m^n:=\{(k,j)\colon k=m,\dots,n;j=1,\dots,2^{k-1}\},
$$
where $m\le n$. 
For $t\in[0,1)$, we also use the branches
$$
  \B(t):=\{(k,j)\colon t\in\Delta_{k-1}^{(j)}\}=\{(k,j)\colon
  \chi_k^{(j)}(t)\not=0 \}.
$$
These are exactly those points of the dyadic tree $\D$ which lie on
one of the infinite paths starting in the root $(1,1)$. 
Figure \ref{fig:dyadic tree} shows a part of the set $\D$, where the
thicker dots stand for the first elements in any of the sets $\B(t)$
for $t\in[\frac3{16},\frac14)$.

We define the following type ideal norms associated with the systems
$\eusm H(\F)$.
\begin{definition}
  For $T\in\frak L(X,Y)$ and for a finite set $\F\subseteq\D$
  denote by $\t(T|\eusm H(\F))$ the least constant $c\ge0$ such that
  \begin{equation}
    \left\| \left. \sum_{\F} Tx_k^{(j)} \chi_k^{(j)} \right|
    L_2 \right\| \le c \left(\sum_{\F} \|x_k^{(j)}\|^2 \right)^{1/2}
  \label{eq3}
  \end{equation}
  for all $\{x_k^{(j)}\colon (k,j)\in\F\}\subseteq X$.

  We call the map
  $$
    \t(\eusm H(\F))\colon T\longrightarrow \t(T|\eusm H(\F))
  $$
  the {\em Haar type ideal norm} associated with the index set $\F$.
  
  We write $\t(X|\eusm H(\F))$ instead of $\t(I_X|\eusm H(\F))$, where 
  $I_X$ is the identity map of the Banach space $X$.
\end{definition}


\section{Comparison of Haar type ideal norms}
\label{sec:3}

In this section we compare Haar type ideal norms associated with 
different index sets $\F$.

The following proposition is an easy consequence of the facts
\eqref{eq2a} and 
(2).
\begin{proposition} \label{prop2}
  Let $n\ge m\ge1$. 
  Then for all $T\in\frak L(X,Y)$ 
  $$
    \t(T|\eusm H(\D_{m+1}^{n+1})) \le \t(T|\eusm H(\D_m^n)).
  $$ 
\end{proposition}

\pf
  We split $\D_{m+1}^{n+1}$ into the two subtrees
  \begin{eqnarray*}
    \S_1 & := & \left\{ (k,j)\in\D_{m+1}^{n+1} \colon 1\le j \le 2^{k-2}
    \right\}, \\
    \S_2 & := & \left\{ (k,j)\in\D_{m+1}^{n+1} \colon 2^{k-2}+1 \le j \le
    2^{k-1} \right\}.
  \end{eqnarray*}
  Note that for $0\le t<1/2$
  $$
    (k,j)\in\S_2\ \text{ implies }\ \chi_k^{(j)}(t)=0 
  $$ 
  and for $1/2\le t<1$
  $$ 
    (k,j)\in\S_1\ \text{ implies }\ \chi_k^{(j)}(t)=0.
  $$
  Hence we may write
  \pagebreak[3]
  \begin{eqnarray*}
    \lefteqn{\int_0^1 \Big\| \sum_{\D_{m+1}^{n+1}} Tx_k^{(j)} \chi_k^{(j)}(t)
    \Big\|^2 dt = } \\
    & = & 
    \int_0^{1/2} \Big\| \sum_{\D_{m+1}^{n+1}}
    Tx_k^{(j)} \chi_k^{(j)}(t) \Big\|^2 dt + \int_{1/2}^1 \Big\|
    \sum_{\D_{m+1}^{n+1}} Tx_k^{(j)} \chi_k^{(j)}(t) \Big\|^2 dt \\
    & = & 
    \int_0^{1/2} \Big\| \sum_{\S_1} Tx_k^{(j)} \chi_k^{(j)}(t)
    \Big\|^2 dt + \int_{1/2}^1 \Big\| \sum_{\S_2} Tx_k^{(j)}
    \chi_k^{(j)}(t) \Big\|^2 dt.
  \end{eqnarray*}
  Substituting $s=2t$ and $h=k-1$ and taking into account formula
  (2), the first integral can be estimated as follows
  \begin{eqnarray*}
    \int_0^{1/2} \Big\| \sum_{\S_1} Tx_k^{(j)} \chi_k^{(j)}(t)
    \Big\|^2 dt 
    & = & 
    \frac12 \int_0^1 \Big\| \sum_{h=m}^n
    \sum_{j=1}^{2^{h-1}} Tx_{h+1}^{(j)} \chi_{h+1}^{(j)} \big(\frac
    s2\big) \Big\|^2 ds \\ & \le & \t(T|\eusm H(\D_m^n))^2 \sum_{h=m}^n
    \sum_{j=1}^{2^{h-1}} \|x_{h+1}^{(j)}\|^2 \\
    & = &
    \t(T|\eusm H(\D_m^n))^2 \sum_{\S_1} \|x_k^{(j)}\|^2.
  \end{eqnarray*}
  Replacing $t$ by $t+\frac12$ and $j$ by $j+2^{k-2}$, and using
  \eqref{eq2a}, the estimate above passes into
  $$
    \int_{1/2}^1 \Big\| \sum_{\S_2} Tx_k^{(j)} \chi_k^{(j)}(t)
    \Big\|^2 dt \le \t(T|\eusm H(\D_m^n))^2 \sum_{\S_2}
    \|x_k^{(j)}\|^2. 
  $$
  which implies the desired inequality.
\qed

By iterated application of Proposition \ref{prop2}, we obtain the
following result.
\begin{corollary}
  Let $m,n\in\N$.
  Then we have for all $T\in\frak L(X,Y)$ that
  $$
    \t(T|\eusm H(\D_{m+1}^{m+n}))\le \t(T|\eusm H(\D_1^n)).
  $$
  \label{lem1}
\end{corollary}

The following transformation and their effects on the Haar functions
will play the crucial role in the proof of Proposition \ref{lem5}.
For $(h,i)\in\D$, we denote by $\phi_h^{(i)}$ the transformation of
$[0,1)$ that interchanges the intervals
$$
  \Delta_{h+1}^{(4i-2)}\ \text{ and }\ \Delta_{h+1}^{(4i-1)}.
$$
More formally
\begin{equation}
  \phi_h^{(i)}(t):=
  \begin{cases}
    t+\frac1{2^{h+1}} & \text{for $t\in\Delta_{h+1}^{(4i-2)}$} \\
    t-\frac1{2^{h+1}} & \text{for $t\in\Delta_{h+1}^{(4i-1)}$} \\
    t & \text{otherwise.}
  \end{cases}
  \label{eq9}
\end{equation}
Lemma \ref{lem2} describes the behavior of the Haar functions under
$\phi_h^{(i)}$. 
What is $\chi_k^{(j)}\circ\phi_h^{(i)}$ for $(k,j)\in\D$\/? 
It turns out that the most interesting cases are those if $(k,j)$
belongs to the fork $F_h^{(i)}:=\{(h,i),(h+1,2i-1),(h+1,2i)\}$.

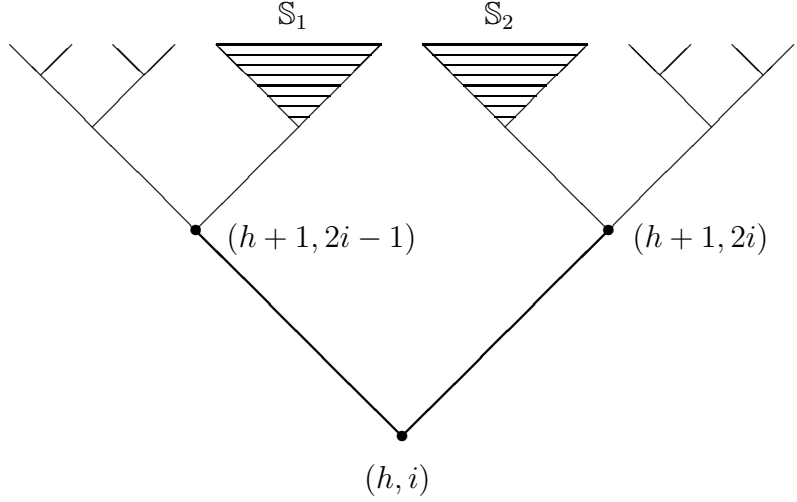
\begin{figure}[h]
  \begin{center}
    \leavevmode
    \setlength{\unitlength}{0.018in}%
    \begin{picture}(228,150)(96,555)
      \thinlines
      \put(270,630){\circle*{3}}    \put(210,570){\circle*{3}}
      \put(150,630){\circle*{3}}
      \put(210,570){\line( 1, 1){114}} \put(300,660){\line(-1, 1){ 24}}
      \put(270,630){\line(-1, 1){ 54}} \put(240,660){\line( 1, 1){ 24}}
      \put(210,570){\line(-1, 1){114}} \put(315,675){\line(-1, 1){  9}}
      \put(285,675){\line( 1, 1){  9}} \put(135,675){\line(-1, 1){  9}}
      \put(105,675){\line( 1, 1){  9}} \put(120,660){\line( 1, 1){ 24}}
      \put(150,630){\line( 1, 1){ 54}} \put(180,660){\line(-1, 1){ 24}}
      \thicklines
      \put(210,570){\line( 1, 1){ 60}} \put(210,570){\line(-1, 1){ 60}}
      \thinlines
      \put(156,684){\line( 1, 0){ 48}} \put(159,681){\line( 1, 0){ 42}}
      \put(162,678){\line( 1, 0){ 36}} \put(165,675){\line( 1, 0){ 30}}
      \put(168,672){\line( 1, 0){ 24}} \put(171,669){\line( 1, 0){ 18}}
      \put(174,666){\line( 1, 0){ 12}} \put(177,663){\line( 1, 0){  6}}
      \put(216,684){\line( 1, 0){ 48}} \put(219,681){\line( 1, 0){ 42}}
      \put(222,678){\line( 1, 0){ 36}} \put(225,675){\line( 1, 0){ 30}}
      \put(228,672){\line( 1, 0){ 24}} \put(231,669){\line( 1, 0){ 18}}
      \put(234,666){\line( 1, 0){ 12}} \put(237,663){\line( 1, 0){  6}}
      \put(199,555){%
        \makebox(0,0)[lb]{$(h,i)$}}
      \put(159,625){%
        \makebox(0,0)[lb]{$(h+1,2i-1)$}}
      \put(277,625){%
        \makebox(0,0)[lb]{$(h+1,2i)$}}
      \put(174,690){%
        \makebox(0,0)[lb]{$\S_1$}}
      \put(234,690){%
        \makebox(0,0)[lb]{$\S_2$}}
    \end{picture}
  \end{center}
  \caption{Part of the dyadic tree $\D_1^n$}
  \label{fig:1}
\end{figure}

If on the other hand the support of $\chi_k^{(j)}$ belongs to one of
the intervals $\Delta_{h+1}^{(4i-2)}$ or $\Delta_{h+1}^{(4i-1)}$ then
obviously $\chi_k^{(j)}\circ\phi_h^{(i)}$ is again a Haar function on
the same level $k$. 
This is exactly the case if $(k,j)$ belongs to one of the subtrees
\begin{eqnarray*}
 \S_1 & := &
 \{(k,j)\colon \Delta_{k-1}^{(j)}\subseteq\Delta_{h+1}^{(4i-2)}\}, \\
 \S_2 & := &
 \{(k,j)\colon \Delta_{k-1}^{(j)}\subseteq\Delta_{h+1}^{(4i-1)}\}.
\end{eqnarray*}
In the remaining cases the Haar functions $\chi_k^{(j)}$ are invariant 
under $\phi_h^{(i)}$. 
Look at Figure \ref{fig:1} to see, how the different index sets
$\F_h^{(i)}$, $\S_1$ and $\S_2$ are related to each other.
\begin{lemma}
  On the fork $\F_h^{(i)}$ the transformation $\phi_h^{(i)}$ acts as
  follows
  \begin{eqnarray*}
    \chi_h^{(i)}\circ \phi_h^{(i)} 
    & = & \big(\chi_{h+1}^{(2i-1)}+\chi_{h+1}^{(2i)}\big)/\sqrt2, \\
    \chi_{h+1}^{(2i-1)}\circ \phi_h^{(i)} 
    & = & \big( \sqrt2\chi_h^{(i)} + \chi_{h+1}^{(2i-1)} - 
    \chi_{h+1}^{(2i)} \big)/2, \\
    \chi_{h+1}^{(2i)}\circ \phi_h^{(i)} 
    & = & \big( \sqrt2\chi_h^{(i)} - \chi_{h+1}^{(2i-1)} +
    \chi_{h+1}^{(2i)} \big)/2.
  \end{eqnarray*}
 Moreover, for $(k,j)\in \S_1$ or $(k,j)\in\S_2$, we get
  $$
    \chi_k^{(j)}\circ \phi_h^{(i)}=
    \begin{cases}
      \chi_k^{(j+2^{k-h-2})} & \text{if $(k,j)\in\S_1$} \\
      \chi_k^{(j-2^{k-h-2})} & \text{if $(k,j)\in\S_2$.}
    \end{cases}
  $$
  The remaining Haar functions $\chi_k^{(j)}$ are invariant under 
  $\phi_h^{(i)}$.
  \label{lem2}
\end{lemma}
\pf
 Looking at Figure \ref{fig:2}, we easily see that
 \begin{eqnarray*}
   \chi_h^{(i)}\circ \phi_h^{(i)} 
   & = & \big( \chi_{h+1}^{(2i-1)}+\chi_{h+1}^{(2i)} \big)/\sqrt2 \\
   (\chi_{h+1}^{(2i-1)} + \chi_{h+1}^{(2i)}) \circ \phi_h^{(i)} 
   & = & \sqrt2 \chi_h^{(i)} \\
   (\chi_{h+1}^{(2i-1)} - \chi_{h+1}^{(2i)}) \circ \phi_h^{(i)} 
   & = & \chi_{h+1}^{(2i-1)}-\chi_{h+1}^{(2i)}.
 \end{eqnarray*}
 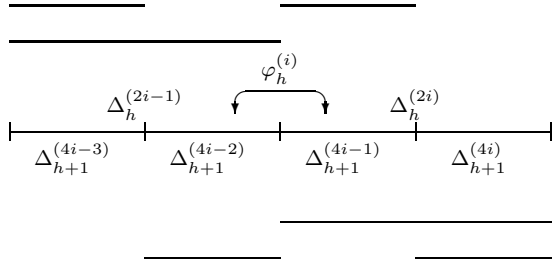
\begin{figure}[h]
   \begin{center}
     \leavevmode
     \unitlength=1.20mm
     \linethickness{0.4pt}
     \begin{picture}(70.00,34.00)(0.00,5.00)
       \put(40.00,23.00){\oval(10.00,3.00)[t]}
       \put(35.00,23.00){\vector(0,-1){1.00}}
       \put(45.00,23.00){\vector(0,-1){1.00}}
       \put(17.00,17.00){%
         \makebox(0,0)[cc]{$\scriptstyle\Delta_{h+1}^{(4i-3)}$}}
       \put(32.00,17.00){%
         \makebox(0,0)[cc]{$\scriptstyle\Delta_{h+1}^{(4i-2)}$}}
       \put(47.00,17.00){%
         \makebox(0,0)[cc]{$\scriptstyle\Delta_{h+1}^{(4i-1)}$}}
       \put(62.00,17.00){\makebox(0,0)[cc]{$\scriptstyle\Delta_{h+1}^{(4i)}$}}
       \put(55.00,23.00){\makebox(0,0)[cc]{$\scriptstyle\Delta_h^{(2i)}$}}
       \put(25.00,23.00){\makebox(0,0)[cc]{$\scriptstyle\Delta_h^{(2i-1)}$}}
       \put(40.00,27.00){\makebox(0,0)[cc]{$\scriptstyle\phi_h^{(i)}$}}
       \put(10.00,30.00){\line(1,0){30.00}}
       \put(10.00,34.00){\line(1,0){15.00}}
       \put(25.00,6.00){\line(1,0){15.00}}
       \put(40.00,10.00){\line(1,0){30.00}}
       \put(55.00,6.00){\line(1,0){15.00}}
       \put(40.00,34.00){\line(1,0){15.00}}
       \put(70.00,20.00){\line(-1,0){60.00}}
       \put(10.00,21.00){\line(0,-1){2.00}}
       \put(25.00,21.00){\line(0,-1){2.00}}
       \put(40.00,21.00){\line(0,-1){2.00}}
       \put(55.00,21.00){\line(0,-1){2.00}}
       \put(70.00,21.00){\line(0,-1){2.00}}
     \end{picture}     
   \end{center}
   \caption{The graphs of the functions 
     $\chi_h^{(i)}$, $\chi_{h+1}^{(2i-1)}$ and $\chi_{h+1}^{(2i)}$}
   \label{fig:2}
 \end{figure}
 
 Solving this system of equations, one gets the assertions for 
 $(k,j)\in\F_h^{(i)}$.

 The other assertions follow easily from the definition of $\phi_h^{(i)}$.
\qed

The following concept turns out to be very useful in the proof of the
main theorem.
\begin{definition}
  Let $\F\subseteq\D$ be a finite subset. 
  The {\em local height} of $\F$ is defined as the maximal number of
  indices in $\F$, that are contained in one branch $\B(t)$. 
  That is, we let
  $$
    \height(\F):=\max_{t\in[0,1)}|\F\cap\B(t)|.
  $$
\end{definition}

The next definition is due to the special behavior of the Haar
functions $\chi_k^{(j)}$ under $\phi_h^{(i)}$ for the indices $(k,j)$
in the fork $\F_h^{(i)}$.
\begin{definition}
  For given $\F\subseteq\D$, we say that an index $(h,i)\in\F$ is 
  {\em $\F$--admissible}, if its successors $(h+1,2i-1)$ and $(h+1,2i)$ do 
  not belong to $\F$.
\end{definition}

For an $\F$--admissible index $(h,i)\in\F$, we assign to $\F$ another
index set $\Phi_h^{(i)}(\F)$ as follows. 
We let
$$
  (h+1,2i-1),(h+1,2i)\in\Phi_h^{(i)}(\F).
$$  
Moreover, we let $(k,j^\ast)\in\Phi_h^{(i)}(\F)$ for all
$(k,j)\in\F\setminus\F_h^{(i)}$, where $j^\ast$ is determined by
$$
  \chi_k^{(j)}\circ\phi_h^{(i)}=\chi_k^{(j^\ast)}.
$$
Note that by Lemma \ref{lem2} such a number $j^\ast$ exists. 

The transformations $\Phi_h^{(i)}$ obviously enjoy the following properties.
\begin{lemma}
 Let $(h,i)\in\F$ be $\F$--admissible. 
 Then
 $$
   |\Phi_h^{(i)}(\F)|=|\F|+1 \text{ and } 
   \height(\Phi_h^{(i)}(\F))=\height(\F).
 $$
\end{lemma}

The next lemma clarifies the importance of the index set
$\Phi_h^{(i)}(\F)$.
\begin {lemma}
  Let $(h,i)\in\F$ be $\F$--admissible. 
  If $f_{\F}$ is given by
  $$ 
    f_{\F}=\sum_{\F}x_k^{(j)}\chi_k^{(j)}\ \text{ with } x_k^{(j)}\in X,
  $$ 
  then there exist elements $y_k^{(j)}\in X$ such that
  $$ 
    f_{\F}\circ\phi_h^{(i)}=\sum_{\Phi_h^{(i)}(\F)}y_k^{(j)}\chi_k^{(j)}.
  $$ 
  Moreover, the $l_2$--sum of the families $(x_k^{(j)})$ and
  $(y_k^{(j)})$ are the same,
  $$
    \sum_{\F}\|x_k^{(j)}\|^2=\sum_{\Phi_h^{(i)}(\F)}\|y_k^{(j)}\|^2.
  $$
  \label{lem4}
\end{lemma}
\pf
 Write $f_{\F}$ in the form
 $$
   f_{\F}=x_h^{(i)}\chi_h^{(i)}+\sum_{\F_0}x_k^{(j)}\chi_k^{(j)}
 $$ 
 with $\F_0=\F\setminus\{(h,i)\}$. 
 In view of Lemma \ref{lem2}
 $$ 
   f_{\F}\circ\phi_h^{(i)} =
   \tfrac1{\sqrt2}x_h^{(i)}\chi_{h+1}^{(2i-1)} +
   \tfrac1{\sqrt2}x_h^{(i)}\chi_{h+1}^{(2i)} +
   \sum_{\F_0}x_k^{(j)}\chi_k^{(j^\ast)}.
 $$ 
 Again, $j^\ast$ is determined by
 $$ 
   \chi_k^{(j)}\circ\phi_h^{(i)}=\chi_k^{(j^\ast)}.
 $$ 
 Note that
 $$  
   \sum_{\F_0}x_k^{(j)}\chi_k^{(j^\ast)} =  
   \sum_{\Phi_h^{(i)}(\F_0)}x_k^{(j^\ast)}\chi_k^{(j)}.
 $$ 
 Hence the elements $y_k^{(j)}$ are given by
 \begin{eqnarray*}
   & y_{h+1}^{(2i-1)}=\tfrac1{\sqrt2}x_h^{(i)},\ 
     y_{h+1}^{(2i)}=\tfrac1{\sqrt2}x_h^{(i)}\ \text{and}\\
   & y_k^{(j)}=x_h^{(j^\ast)}\ \text{ for }\ (k,j)\in\Phi_h^{(i)}(\F_0),
 \end{eqnarray*}
 \nopagebreak[4]
 and their $l_2$--sum can easily be computed.
\qed
\pagebreak[3]

Note that an index set $\F$ of local height $n$ can be arbitrarily
large.  
However, the following lemma shows that, using the transforms
$\Phi_h^{(i)}$, such an index set can be compressed so as to fit into
a certain dyadic tree $\D_{m+1}^{m+n}$. 
The idea is to replace an admissible index $(h,i)\in\F$ by its
successors $(h+1,2i-1)$ and $(h+1,2i)$.
\begin{lemma}
 Let $n:=\height(\F)$. 
 Then there exist $m\in\N$ and a finite sequence of transforms
 $$ 
   \F=:\F_0 \stackrel{\Phi_{h_1}^{(i_1)}}{\longrightarrow}
   \F_1 \stackrel{\Phi_{h_2}^{(i_2)}}{\longrightarrow}
   \dots \stackrel{\Phi_{h_N}^{(i_N)}}{\longrightarrow}
   \F_N
 $$ 
 such that
 $$ 
   \F_N\subseteq\D_{m+1}^{m+n}.
 $$
 \label{lem3}
\end{lemma}
\pf
  Since $\F$ is finite, we can choose $m\in\N$ with
  $\F\subseteq\D_1^{m+n}$. 
  As often as possible, we successively apply transforms
  $\Phi_{h_l}^{(i_l)}$ for $\F_{l-1}$--admissible indices $(h_l,i_l)$ 
  $$
    \F_{l-1} \stackrel{\Phi_{h_l}^{(i_l)}}{\longrightarrow} \F_l.
  $$ 
  In order to guarantee that $\F_l$ is still contained in
  $\D_{m+1}^{m+n}$, we always assume that $h_l<m+n$.
 
  Since
  $$ 
    |\F_l|=|\F_{l-1}|+1\ \text{ and }\ |\F_l|\le|\D_1^{m+n}|,
  $$
  this process terminates after a finite number of steps. 
  Hence, the last index set $\F_N$ does not contain any
  $\F_N$--admissible index $(h,i)$ with $h<m+n$. 
 
  Choose $h$ such that
  $$ 
    \F_N\subseteq\D_h^{m+n}\ \text{ and }\ \F_N\not\subseteq\D_{h+1}^{m+n}.
  $$
  \pagebreak[3]
  If $h=m+n$, then
  $\F_N\subseteq\D_{m+n}^{m+n}\subseteq\D_{m+1}^{m+n}$. Thus we are
  done. 
  Otherwise, take any index $(h,j_0)\in\F_N$ on the lowest
  level. 
  Since $h<m+n$ and $(h,j_0)$ is not $\F_N$--admissible, at least one
  of its successors, say $(h+1,j_1)$, must belong to $\F_N$.
  In this way, we find a sequence
  $$ 
    (h,j_0),\ (h+1,j_1),\dots,(m+n,j_{m+n-h})\in\F_N
  $$ 
  of length $m+n-h+1$, that belongs to some branch $\B(t)$. 
  Hence
  $$ 
    m+n-h+1\le n=\height(\F).
  $$
  Finally, we conclude from $m+1\le h$ that
  $$ 
    \F_N\subseteq\D_h^{m+n}\subseteq\D_{m+1}^{m+n}.
  $$
\qed

\pagebreak[3]
We can now formulate the most interesting result of this section.
\begin{proposition}
  Let $\F\subseteq\D$ be a finite subset of local height $n$. 
  Then there exists $m\in\N$ such that for all $T\in\frak L(X,Y)$
  $$ 
    \t(T|\eusm H(\F)) \le \t(T|\eusm H(\D_{m+1}^{m+n})).
  $$
  \label{lem5}
\end{proposition}
\pf
  Consider the transforms
  $\Phi_{h_1}^{(i_1)},\dots,\Phi_{h_N}^{(i_N)}$ constructed in the
  previous Lemma. 
  Let $\phi_{h_1}^{(i_1)},\dots,\phi_{h_N}^{(i_N)}$ denote the
  associated bijections defined in \eqref{eq9}.
 
 Given
 $$ 
   f_{\F}=\sum_{\F}x_k^{(j)}\chi_k^{(j)},
 $$
 we let
 $$ 
   f_{\F_l}:=f_{\F_{l-1}}\circ\phi_{h_l}^{(i_l)}\text{ and } 
   f_{\F_0}:=f_{\F}.
 $$ 
 Write
 $$ 
   f_{\F_l}=\sum_{\F_l}x_k^{(j,l)}\chi_k^{(j)}.
 $$
 Since the transformations $\phi_h^{(i)}$ are measure preserving, we have
 $$
   \Big\|\sum_{\F}Tx_k^{(j)}\chi_k^{(j)}\Big|L_2\Big\| \! = \!
   \|f_{\F\unskip_0}|L_2\| \! = \dots = \! \|f_{\F\unskip_N}|L_2\| \! = \!
   \Big\|\sum_{\F\unskip_N}Tx_k^{(j,N)}\chi_k^{(j)}\Big|L_2\Big\|.
 $$ 
 Moreover, Lemma \ref{lem4} yields that
 $$
   \sum_{\F}\|x_k^{(j)}\|^2 = \sum_{\F_0}\|x_k^{(j,0)}\|^2 = \dots =
   \sum_{\F_N}\|x_k^{(j,N)}\|^2.
 $$ 
 Since $\F_N\subseteq\D_{m+1}^{m+n}$, we have
 $$ 
   \Big\|\sum_{\F_N}Tx_k^{(j,N)}\chi_k^{(j)}\Big|L_2\Big\| \le 
   \t(T|\eusm H(\D_{m+1}^{m+n}))
   \left(\sum_{\F_N}\|x_k^{(j,N)}\|^2\right)^{1/2}.
 $$ 
 Hence
 $$ 
   \Big\|\sum_{\F}Tx_k^{(j)}\chi_k^{(j)}\Big|L_2\Big\| \le 
   \t(T|\eusm H(\D_{m+1}^{m+n}))
   \left(\sum_{\F}\|x_k^{(j)}\|^2\right)^{1/2}.
 $$
\qed

We now summarize the results of this section.
\begin{theorem}
  Let $\F\subseteq\D$ be a finite subset of local height $n$. 
  Then we have for all $T\in\frak L(X,Y)$ that
  $$ 
    \t(T|\eusm H(\F)) \le \t(T|\eusm H(\D_1^n)).
  $$
  \label{th:comparision}
\end{theorem}

\pf
  The assertion follows immediately from Corollary \ref{lem1} and 
  Proposition \ref{lem5}.
\qed

\begin{remark}
  Using the same methods as in the proof of Proposition \ref{lem5}, one can 
  show that
  $$
    \t(T|\eusm H(\F))=\t(T|\eusm H(\D_1^n))
  $$
  for all $T\in\frak L(X,Y)$, if $\F$ has {\em exact local height} $n$, 
  i.e.~if
  $$
    |\F\cap\B(t)|=n
  $$
  for all $t\in[0,1)$.
\end{remark}

As a consequence of the remark above, we get the following result.

\begin{corollary}
  Let $m,n\in\N$.
  Then we have for all $T\in\frak L(X,Y)$ that
  $$
    \t(T|\eusm H(\D_{m+1}^{m+n})) = \t(T|\eusm H(\D_1^n)).
  $$
\end{corollary}


\section{A combinatorial lemma}
\label{sec:3a}

In this section we provide a lemma, which is needed in the proof of
Theorem \ref{th:local}.
\begin{lemma}
  Let $\F\subseteq\D_1^n$ be an index set of local height $l\le n$. 
  If $|\F|<2^l-1$ then there exists an index
  $(k,j)\in\D_1^n\setminus\F$ such that
  $$
    \height(\F\cup\{(k,j)\}) = l.
  $$
\end{lemma}
\pf
  For $l=1$ the assertion is trivial. 
  Now assume that for $l-1$ the lemma is true. 
  To prove the lemma for $l$, we use induction over $n\ge l$.

  If $n=l$ then $\D_1^n$ has exactly $2^n-1=2^l-1$ elements and each
  subset of $\D_1^n$ has local height less than or equal to $l$.
  Hence, we may take any element $(k,j)\in\D_1^n\setminus\F$. 
  Since $|\F|<2^l-1$, the latter set is certainly nonempty.

  Now assume that $\F\subseteq\D_1^n$ and $n>l$. 
  Note that the subtrees
  \begin{eqnarray*}
    \S_1 
    & := & 
    \left\{(k,j)\in\D_2^n\colon 1\le j\le2^{k-2}\right\},\\
    \S_2 
    & := & 
    \left\{(k,j)\in\D_2^n\colon 2^{k-2}+1\le j\le2^{k-1}\right\}
  \end{eqnarray*}
  can be canonically identified with $\D_1^{n-1}$. 
  Let
  $$
    \F_1:=\F\cap\S_1\ \text{ and }\ \F_2:=\F\cap\S_2
  $$
  and consider them as subsets of $\D_1^{n-1}$ via the identification
  above.

  If $(1,1)\notin\F$ then
  $$
    \height(\F_1)\le l\ \text{ and }\ |\F_1|\le|\F|<2^l-1.
  $$
  By the induction hypothesis, there exists
  $(k,j)\in\S_1\setminus\F_1$ such that
  $$
    \height(\F_1\cup\{(k,j)\})\le l. 
  $$
  Since $(1,1)\notin\F$, this implies that also
  $$ 
    \height(\F\cup\{(k,j)\})\le l.
  $$
  
  If however $(1,1)\in\F$ then 
  $$ 
    \height(\F_1)\le l-1\ \text{ and }\ \height(\F_2)\le l-1.
  $$
  Moreover, without loss of generality, we may assume that
  $$ 
    |\F_1|\le\frac{|\F-1|}2<2^{l-1}-1. 
  $$
  Since we already know that the lemma is true for $l-1$, there exists
  $(k,j)\in\D_1^{n-1}\setminus\F_1$ such that
  $$ 
    \height(\F_1\cup\{(k,j)\})\le l-1. 
  $$
  Since $(1,1)\in\F$, this implies that
  $$ 
    \height(\F\cup\{(k,j)\})\le l.
  $$
\qed

By repeated application of the previous lemma, one easily checks the
following statement.
\begin{corollary} \label{cor:combinatoric}
  Let $\F\subseteq\D_1^n$ be an index set satisfying $\height(\F)\le
  l\le n$. 
  If $|\F|<2^l-1$ then there exists an index set
  $\F_0\subseteq\D_1^n\setminus\F$ of cardinality $2^l-1-|\F|$ such
  that
  $$
    \height(\F\cup\F_0) \le l.
  $$
\end{corollary}


\section{Haar type and weak Haar type}
\label{sec:4}

First of all, let us introduce variants of the ideal norms $\t(\eusm
H(\D_1^n))$ defined in Section~\ref{sec:2}.

\begin{definition}
  For $T\in\frak L(X,Y)$ and for $1\le p\le2$ denote by $\t_p(T|\eusm
  H(\D_1^n))$ the least constant $c\ge0$ such that
  $$
    \bigg\| \sum_{\D_1^n} Tx_k^{(j)} \chi_k^{(j)} \bigg| L_p \bigg\|
    \le c \left( \sum_{k=1}^n \bigg\| \sum_{j=1}^{2^{k-1}} x_k^{(j)}
    \chi_k^{(j)} \bigg| L_p \bigg\|^p \right)^{1/p} 
  $$ 
  for all $(x_k^{(j)}) \subseteq X$.
\end{definition}

Note that
$$ 
  \t_2(\eusm H(\D_1^n)) = \t(\eusm H(\D_1^n)). 
$$

The following theorem allows to easily determine the asymptotic
behavior of the ideal norms $\t_p(\eusm H(\D_1^n))$. 
It is a refinement of a theorem of Pisier in \cite{pis75} and is due
to Gei\ss{} \cite{gei}.
\begin{theorem}
  Let $n\in\N$ and $1\le p\le2$ be fixed. 
  For $T\in\frak L(X,Y)$ assume that
  $$ 
    \bigg\| \sum_{\D_1^n} Tx_k^{(j)} \chi_k^{(j)} \bigg| L_1 \bigg\|
    \le c_n \, \left\| \left. \left( \sum_{k=1}^n \bigg\|
    \sum_{j=1}^{2^{k-1}} x_k^{(j)} \chi_k^{(j)} \bigg\|^p \right) ^{1/p}
    \right| L_\infty \right\|
  $$ 
  for all $(x_k^{(j)}) \subseteq X$.
  Then it follows that
  $$ 
    \t_p(T|\eusm H(\D_1^n))\le c\,c_n.
  $$ 
  Here $c$ is a universal constant not depending on $n$ nor $p$.
  \label{th3}
\end{theorem}

The next definitions are motivated by the corresponding concepts in
the case of Rademacher functions.
\begin{definition}
  For $1\leq p\leq2$, we say that an operator $T\in\frak L(X,Y)$ is of
  {\em Haar type\/} $p$, if there exists a constant $c\geq0$ such that
  $$
    \t_p(T|\eusm H(\D_1^n)) \leq c 
  $$ 
  for all $n\in\N$.

  For $1\leq p<2$, we say that an operator $T\in\frak L(X,Y)$ is of
  {\em weak Haar type\/} $p$, if there exists a constant $c\geq0$ such
  that 
  $$ 
    \t(T|\eusm H(\D_1^n)) \leq c\,n^{1/p-1/2} 
  $$ 
  for all $n\in\N$.
\end{definition}

\pagebreak[3]
\def\remarkname{Remarks}
\begin{remark}
  \begin{itemize}
  \item 
    Theorem \ref{th3} ensures that for $p<2$ Haar type $p$ implies
    weak Haar type $p$. 
    Indeed, note that 
    \begin{eqnarray*}
      \lefteqn{\bigg\| \sum_{\D\unskip_1^n} Tx_k^{(j)}\chi_k^{(j)}
        \bigg| L_1 \bigg\|
      \le 
        \t_p(T|\eusm H(\D_1^n))\, n^{1/p-1/2}\, \times} \\
      & & \left\|\left. \bigg(\sum_{k=1}^n \bigg\| \sum_{j=1}^{2^{k-1}}
      x_k^{(j)} \chi_k^{(j)} \bigg\| ^2 \bigg)^{1/2} \right| L_\infty
      \right\|
    \end{eqnarray*}
    and hence by Theorem \ref{th3} 
    $$ 
      \t(T|\eusm H(\D_1^n)) \le c\, \t_p(T|\eusm H(\D_1^n))\,
      n^{1/p-1/2}.
    $$ 
  \item 
    In the case $p=2$ the corresponding definition of weak Haar
    type $2$ would obviously coincide with that of Haar type
    $2$. That's why, we consider weak Haar type $p$ only for $p<2$. 
  \item 
    In Pisier's work \cite{pis75} it was shown that a Banach space
    $X$ is of Haar type $p$ exactly if it admits an equivalent
    $p$--smooth renorming. A similar statement also holds for
    operators; see also \cite{BEA85} or the forthcoming book
    \cite{PIW} for this connection.
  \item 
    Pisier in \cite{pis86} also introduced a concept of martingale
    type $p$, which again is equivalent to Haar type $p$. This follows
    from the considerations in \cite{pis75}.
  \end{itemize}
\end{remark}

We can now prove the main theorem of this article.

\begin{theorem}\label{th:main}
  If an operator $T\in\frak L(X,Y)$ is of weak Haar type $p$ for some
  $1\leq p<2$ then it is of Haar type $r$ for all $r<p$.
\end{theorem}

\pf
  The proof follows essentially that of Sublemma 3.1. in Pisier
  \cite{pis75}.  
  The main new ingredient is the use of the results of Section
  \ref{sec:3}.

  Let $r<p$ and $(x_h^{(i)})\subseteq X$. 
  Set 
  $$ 
    S_r:= \Bigg\| \left(\sum_{k=1}^n \bigg\| \sum_{j=1}^{2^{k-1}}
    x_k^{(j)} \chi_k^{(j)} \bigg\|^r \right)^{1/r} \Bigg| L_\infty
    \Bigg\|
  $$ 
  and for $l=1,2,\dots$ define
  $$
    \F_l:=\left\{ (k,j)\in\D_1^n \colon \frac {S_r}{2^{l/r}} < 2^{(k-1)/2} 
    \|x_k^{(j)}\| \le \frac {S_r}{2^{(l-1)/r}} \right\}.
  $$ 
  Then it follows that
  \begin{equation}
    \label{eq:partition}
    \sum_{\D_1^n} x_k^{(j)} \chi_k^{(j)} = \sum_{l=1}^\infty \sum_{\F_l} 
    x_k^{(j)} \chi_k^{(j)}.
  \end{equation}
  Moreover the sets $\F_l$ have local height less than $2^l$. 
  To see this last fact, choose $t\in[0,1)$ and note that
  \begin{eqnarray*}
    S_r^r & \ge & \sum_{\D_1^n} \|x_k^{(j)} \chi_k^{(j)}(t) \|^r \\ 
    & \ge & \sum_{\F_l\cap \B(t)} \| x_k^{(j)} \chi_k^{(j)}(t) \|^r >
    \sum_{\F_l\cap \B(t)} \left( \frac {S_r}{2^{l/r}} \right)^r =
    |\F_l\cap \B(t)|\, S_r^r \, 2^{-l}.
  \end{eqnarray*}
  This shows that
  \begin{equation}
    |\F_l\cap \B(t)|<2^l.
    \label{eq:local_height}
  \end{equation}

  By Theorem \ref{th:comparision} and the definition of weak Haar type
  $p$, there exists a constant $c$ such that
  \begin{equation}
    \t(T|\eusm H(\F_l)) \le \t(T|\eusm H(\D_1^{2^l})) \le
    c\,2^{l(1/p-1/2)}.
    \label{eq:est_F_l}
  \end{equation}
  
  Since $r<p$ it follows that 
  $$ 
    \sum_{l=1}^\infty 2^{l(1/p-1/r)} = c_{pr} < \infty. 
  $$ 
  Hence \eqref{eq:est_F_l} implies
  \begin{equation}
    \sum_{l=1}^\infty \t(T|\eusm H(\F_l))\, 2^{l(1/2-1/r)} \le c \,
    c_{pr}.
    \label{eq:powersum_finite}
  \end{equation}
  
  Moreover we have
  \begin{eqnarray}
    \sum_{\F_l} \|x_k^{(j)}\|^2
    & = &
    \sum_{\F_l} \int_0^1 \|x_k^{(j)}\chi_k^{(j)}(t)\|^2 dt
    =
    \int_0^1
    \sum_{\F_l\cap\B(t)} \|x_k^{(j)} \chi_k^{(j)}(t) \|^2 dt
    \nonumber\\
    & = &
    \int_0^1\sum_{\F_l\cap\B(t)} 2^{k-1}\|x_k^{(j)}\|^2 dt \\
    & \le &
    \int_0^1 |\F_l\cap\B(t)| \left( \tfrac {S_r}{2^{(l-1)/r}}
    \right)^2 dt
    \le
    2^{2/r} 2^{l(1-2/r)}S_r^2.\nonumber
    \label{eq:estimat_sum_squares}
  \end{eqnarray}
  
  \pagebreak[3]
  Now \eqref{eq:partition}, the definition of $\t(T|\eusm H(\F_l))$,
  \eqref{eq:estimat_sum_squares}, and \eqref{eq:powersum_finite} imply
  \begin{eqnarray*}
    \Big\| \sum_{\D_1^n} Tx_k^{(j)} \chi_k^{(j)} \Big| L_2 \Big\| 
    & \le & 
    \sum_{l=1}^\infty \Big\| \sum_{\F_l} Tx_k^{(j)} \chi_k^{(j)}
    \Big| L_2 \Big\| \\
    & \le &
    \sum_{l=1}^\infty \t(T|\eusm H(\F_l)) \left(
    \sum_{\F_l} \|x_k^{(j)} \|^2 \right)^{1/2} \\ 
    & \le & 2^{1/r}
    \sum_{l=1}^\infty \t(T|\eusm H(\F\unskip_l)) 2^{l(1/2-1/r)} S_r \le
    2^{1/r} c \, c_{pr} \, S_r.
  \end{eqnarray*}
  Finally Theorem \ref{th3} gives
  $$ 
    \t_r(T|\eusm H(\D_1^n)) \le 2^{1/r}\,c \, c_{pr}, 
  $$ 
  which completes the proof.
\qed

Looking at this proof more closely and exploiting Corollary
\ref{cor:combinatoric}, we can even prove the following local variant
of the previous result:
\begin{theorem}\label{th:local}
  If an operator $T\in\frak L(X,Y)$ is of weak Haar type $p$ for some
  $1\leq p<2$ then there exists a constant $c\geq0$ such that 
  $$
    \t_p(T|\eusm H(\D_1^n)) \leq c\, (1+\log n) 
  $$ 
  for all $n\in\N$.
\end{theorem}
\pf
  Let $(x_k^{(j)})\subseteq X$. 
  As in the previous proof define
  $$ 
    \F_l:=\left\{ (k,j)\in\D_1^n \colon \frac {S_p}{2^{l/p}} < 2^{(k-1)/2}
    \|x_k^{(j)}\| \le \frac {S_p}{2^{(l-1)/p}} \right\}. 
  $$ 
  Again \eqref{eq:partition} and \eqref{eq:local_height} hold. 
  Let $m$ be such that $2^m \leq n<2^{m+1}$. 
  We now use Corollary \ref{cor:combinatoric} to construct
  modifications $\F_1',\dots,\F_{m+1}'$ of the sets $\F_l$ of maximal
  cardinality.
  This is done inductively.  
  If $|\F_1|\ge3$ then we let $\F_1':=\F_1$.
  Otherwise, by Corollary~\ref{cor:combinatoric}, there exist a subset
  $\F_1''\subseteq\D_1^n\setminus\F_1$ of cardinality $3-|\F_1|$ such
  that $\height(\F_1\cup\F_1'')\le2$.
  Defining $\F_1':=\F_1\cup\F_1''$, we get that
  $$ 
    \height(\F_1')\le2^1\ \text{ and }\ |\F_1'| \ge 2^{2^1}-1. 
  $$
  Now assume that $l\le m$ and that $\F_1',\dots,\F_{l-1}'$ are
  already defined such that
  $$
    \begin{array}{rlrl}
      \height(\F_{l-1}')       & \le 2^{l-1}, & 
      \sum_{k=1}^{l-1} |\F_k'| & \ge 2^{2^{l-1}}-1, \\
      \F_{l-1}'                & \subseteq \D_1^n \setminus
      \bigcup_{k=1}^{l-2} \F_k, & 
      \bigcup_{k=1}^{l-1} \F_k & \subseteq \bigcup_{k=1}^{l-1} \F_k'.
    \end{array}
  $$
  If
  $$ 
    \sum_{k=1}^{l-1}|\F_k'| + |\F_l| \ge 2^{2^l}-1 
  $$ 
  then we let $\F_l':=\F_l$.
  Otherwise, by Corollary \ref{cor:combinatoric} (note that $2^l\le
  2^m \le n$), there exists a subset
  $\F_l''\subseteq\D_1^n\setminus\F_l$ of cardinality
  $2^{2^l}-1-|\F_l|$ such that $\height(\F_l\cup\F_l'')\le2^l$.
  Defining 
  $$
    \F_l':=(\F_l\cup\F_l'') \setminus \bigcup_{k=1}^{l-1}\F_k',
  $$
  we get that $\height(\F_l')\le2^l$, and
  $$ 
    \sum_{k=1}^l |\F_k'| \ge \sum_{k=1}^{l-1} |\F_k'| +
    (|\F_l|+|\F_l''|) - \sum_{k=1}^{l-1}|\F_k'| = 2^{2^l}-1. 
  $$
  Moreover, by construction
  $$ 
    \F_l'\subseteq \D_1^n \setminus \bigcup_{k=1}^{l-1} \F_k
    \ \text{ and }\ 
    \bigcup_{k=1}^l\F_k \subseteq \bigcup_{k=1}^l \F_k'.
  $$
  The last condition ensures that for $(k,j)\in\F_l'$ we have
  $(k,j)\notin \bigcup_{h=1}^{l-1}\F_h$ and hence
  $$ 
    2^{(k-1)/2}\|x_k^{(j)}\|\le \frac {S_p}{2^{(l-1)/p}}. 
  $$ 
  This construction yields a sequence $\F_1',\dots,\F_m'$ of length
  $m$. 
  Finally, we let 
  $$
    \F_{m+1}':=\D_1^n\setminus\bigcup_{l=1}^m \F_l'.
  $$ 
  Obviously $\height(\F_{m+1}') \le n < 2^{m+1}$. 
  Moreover, for $(k,j)\in\F_{m+1}'$,
  $$ 
    2^{(k-1)/2}\|x_k^{(j)}\|\le \frac {S_p}{2^{m/p}}.
  $$
    
  Since again $\height(\F_l')<2^l$ for $l=1,\dots,m+1$, we have
  \eqref{eq:est_F_l} for $\F_l'$.
  Moreover
  $$
    \sum_{l=1}^{m+1} \t(T|\eusm H(\F_l'))\,2^{l(1/2-1/p)} \le c\,
    (m+1) \le c\, (1+\log n).
  $$
  Furthermore, Inequality \eqref{eq:estimat_sum_squares} holds
  with $r$ replaced by $p$. 
  Now, since 
  $$ 
    \bigcup_{l=1}^{m+1} \F_l' = \D_1^n, 
  $$ 
  we get
  \begin{eqnarray*}
    \Big\| \sum_{\D_1^n} Tx_k^{(j)} \chi_k^{(j)} \Big| L_2 \Big\| 
    & \le &
    \sum_{l=1}^{m+1} \Big\| \sum_{\F_l'} Tx_k^{(j)} \chi_k^{(j)}
    \Big| L_2 \Big\|\\
    & \le &
    \sum_{l=1}^{m+1} \t(T|\eusm H(\F_l')) \left(
    \sum_{\F_l'} \|x_k^{(j)} \|^2 \right)^{1/2} \\ 
    & \le &
    2^{1/p} \sum_{l=1}^{m+1} \t(T|\eusm H(\F_l')) 2^{l(1/2-1/p)} S_p\\
    & \le &
    2^{1/p} \, c\,(1+\log n)\, S_p.
  \end{eqnarray*}
  A glance at Theorem \ref{th3} completes the proof.
\qed


\section{Final remarks}
\label{sec:final}

To show that the notions of Haar type $p$ and weak Haar type $p$ do
not coincide, we provide the following example.

Let $1<p<2$ and $\sigma_k:=k^{-1/p'}$.
Consider the diagonal operator $D_s\colon l_1\to l_1$ associated with
the sequence $s=(\sigma_k)$, which is given by 
$$
  x=(\xi_k) \longmapsto D_sx=(\sigma_k\xi_k).
$$
Then
$$
  \t(D_s|\eusm H(\D_1^n)) = \left( \sum_{k=1}^n k^{-2/p'}
  \right)^{1/2} \le n^{1/p-1/2}
$$
and hence $D_s$ is of weak Haar type $p$.
On the other hand
$$
  \t_p(D_s|\eusm H(\D_1^n)) = \left( \sum_{k=1}^n k^{-1}
  \right)^{1/p'} \ge \frac 12(1+\log n)^{1/p'}
$$
and hence $D_s$ is not of Haar type $p$.
The tedious calculations are carried out in \cite{wen94}.

However, for {\em spaces\/} the situation seems to be more complicated.
\begin{problem}
  Given $1<p<2$, does there exist a Banach space $X$ that is of weak
  Haar type $p$ but not of Haar type $p$\/?
\end{problem}
In the Rademacher case, such examples are given by the well known
modifications of the Tsirelson space (see Tzafriri \cite{tza79}) which
are even Banach lattices.
It can also be shown, that for Banach lattices the concepts of
Rademacher and Haar type are the same.
However, it seems to be unknown, whether the same is true for the
corresponding weak properties.

The examples above also suggest that a stronger version of Theorem
\ref{th:main} is true.
\begin{problem}
  If $T$ is of weak Haar type $p$, does there exist a constant
  $c\ge0$ such that
  $$
    \t_p(T|\eusm H(\D_1^n)) \le c\, (1+\log n)^{1/p'}\/?
  $$
\end{problem}

One more problem remains open. 
In the introduction it was mentioned that martingale type $p$ and Haar
type $p$ are equivalent properties.
\begin{problem}
  Is it true that weak martingale type $p$ and weak Haar type $p$ are
  the same?
\end{problem}


{\sc Mathematical Institute, FSU Jena, 07740 Jena, Germany}

{\sl E--mail address}: {\tt wenzel@minet.uni-jena.de}

\end{document}